# An Explicit Formula For The Divisor Function

N. A. Carella

**Abstract**: The details for the construction of an explicit formula for the divisors function $d(n) = \#\{ d \mid n \}$ are formalized in this article. This formula facilitates a unified approach to the investigation of the error terms of the divisor problem and circle problem.

**Mathematics Subject Classification:** 11N37, 11N56, 11P21, 11M06.

**Keywords:** Divisor Function, Asymptotic Formula, Divisor Problem, Circle Problem.

## 1. Introduction

Let $n \in \mathbb{N} = \{ 0, 1, 2, 3, \ldots \}$ be an integer, and let $d(n) = \#\{ d : d \text{ divides } n \}$ be the number of divisors counting function. The average order $\sum_{n \leq x} d(n)$ of the divisors counting function is given by

$$\sum_{n \leq x} d(n) = (\log x + 2\gamma - 1)x + O(x^{1/2}), \tag{1}$$

where $\gamma = .5772\ldots$ is Euler constant. The summatory function $\sum_{n \leq x} d(n)$ has the weighted Voronoi formula

$$\frac{x^{v-1}}{\Gamma(v)} \sum_{n \leq x} \left(1 - \frac{n}{x}\right)^{v-1} d(n) = \frac{x^{v-1}}{4\Gamma(v)} + \frac{x^v}{\Gamma(v+1)}\left(\log x + \gamma - \psi(1+v)\right) - 2\pi x^v \sum_{n \geq 1} d(n) F_v(4\pi\sqrt{nx}), \tag{2}$$

where $v \geq 1$ is a parameter, $F_v(z)$ is a linear combination of Bessel functions such as $K_v(z)$ and $Y_v(z)$, $\psi(z) = \Gamma'(z)/\Gamma(z) = -\gamma + \sum_{n \leq z-1} n^{-1}$ is the digamma function, and $-\psi(1) = \gamma$, see [IV, p. 89], [KC, p. 194]. The most common is the special case $v = 1$, known as the Voronoi formula

$$\sum_{n \leq x} d(n) = \frac{1}{4} + \left(\log x + 2\gamma - 1\right)x - \frac{2x^{1/2}}{\pi}\sum_{n \geq 1}\frac{d(n)}{n^{1/2}}\left(K_1(4\pi(nx)^{1/2}) + 2^{-1}\pi Y_1(4\pi(nx)^{1/2})\right) \tag{3}$$

for a number $x \in \mathbb{R} - \mathbb{N}$, confer [IV, p. 83], [KC, p. 195]. Most of the research and investigations on the divisors summatory function $\sum_{n \le x} d(n)$ are based on these formulas.

Evidently, the Voronoi formula is not expressed in terms of the zeros of the zeta function. The reason being that the generating function $\sum_{n \ge 1} d(n)n^{-s} = \zeta(s)^2$ does not have the necessary analytical structures to facilitate the construction of an explicit formula of the divisor summatory function $\sum_{n \le x} d(n)$.

The proof of the main result is on the next Subsection. Section 2 has the materials for a few generating functions. Section 3 has the calculations of the explicit formula for the finite sum $\sum_{n \le x} 2^{\omega(n)}$. The other Sections are references, which include various other related results.

## 1.1 Explicit Formula

An indirect approach of constructing such formula is formulated here. This apparently new explicit formula expresses the arithmetic nature of the summatory function $\sum_{n \le x} d(n)$ in terms of the zeros of the zeta function $\zeta(s)$.

**_Theorem_ 1.1.** For a large number $x \in \mathbb{R} - \mathbb{N}$, the average order of the divisors counting function has the explicit formula

$$\sum_{n \le x} d(n) = -\frac{\pi^2}{12} + \left(\log x + 2\gamma - 1\right)x + \frac{\pi^2}{3}\sum_{\rho}\frac{\zeta(\rho/2)^2}{\rho\zeta'(\rho)}x^{\rho/2} - \frac{\pi^2}{6}\sum_{n \ge 0}\frac{\zeta(-2n-1)^2}{(2n+1)\zeta'(-2n-1)}x^{-2n-1},$$

(4)

where $\rho$ ranges over the nontrivial zeros of the zeta function $\zeta(s)$.

**_Proof_**: Rewrite the generating function $\sum_{n \ge 1} 2^{\omega(n)}n^{-s} = \zeta(s)^2/\zeta(2s)$, see Lemma 2.1, of the squarefree divisors counting function $\#\{d \,|\, n : \mu(n) \ne 0\} = 2^{\omega(n)}$ in the form

$$\sum_{n \ge 1}\frac{d(n)}{n^s} = \zeta(2s)\sum_{n \ge 1}\frac{2^{\omega(n)}}{n^s}.$$

(5)

Setting $s = 1$, this immediately implies that

$$\sum_{n \le x} d(n) = \zeta(2)\sum_{n \le x} 2^{\omega(n)} + E(x),$$

(6)

where $E(x) = O(x\log x)$ is a correction factor. Now, let $x \in \mathbb{R} - \mathbb{N}$, and apply Theorem 3.1, and Lemma 7.3 to get the expression

$$\sum_{n \le x} d(n) = \left( \log x + 2\gamma - 1 \right) x + O(x^{1/2})$$

$$= \zeta(2) \sum_{n \le x} 2^{\omega(n)} + E(x)$$

$$= -\frac{\pi^2}{12} + \left( \log x + 2\gamma - 1 - 12\zeta'(2)/\pi^2 \right) x$$

$$+ \frac{\pi^2}{3} \sum_{\rho} \frac{\zeta(\rho/2)^2}{\rho\zeta'(\rho)} x^{\rho/2} + \frac{\pi^2}{6} \sum_{n \ge 0} \frac{\zeta(-2n-1)^2}{(2n+1)\zeta'(-2n-1)} x^{-2n-1} + E(x) .$$

(7)

Since both sides most have the same asymptotic order (main term), it is clear that

$$\sum_{n \le x} d(n) \sim \left( \log x + 2\gamma - 1 \right) x$$

$$\sim \zeta(2) \sum_{n \le x} 2^{\omega(n)} + E(x)$$

$$\sim \left( \log x + 2\gamma - 1 - 12\zeta'(2)/\pi^2 \right) x + E(x) .$$

(8)

Thus, it follows that the correction factor $E(x) = (12\zeta'(2)/\pi^2) x$. ∎

The earlier formulas as Lemma 4.1 and 7.3 for $\sum_{n \le x} d(n)$ and $\sum_{n \le x} 2^{\omega(n)}$ respectively, also shows the relationship $\sum_{n \le x} d(n) = (12\zeta'(2)/\pi^2)x + \zeta(2)\sum_{n \le x} 2^{\omega(n)}$ between these two arithmetic functions.

Another way of verifying the asymptotic formula (8) involves the calculations of the derivatives of both asymptotic formulas in (7). This concludes with $E(x) = (12\zeta'(2)/\pi^2)x + c$, $c > 0$ constant.

The explicit formula $D_0(x)$ in Theorem 1.1 is a harmonic expansion of the average order $D(x) = \sum_{n \le x} d(n)$ of the divisors counting function $d(n) = \sum_{d \mid n} 1$. Basically, the explicit formula $D_0(x)$ approximates the discontinued function $D(x)$ by an infinite sum of harmonic functions $\left\{ x^{\rho/2} = x^{\sigma/2} e^{it\log x/2} : \zeta(\rho) = 0 \right\}$. Since the summatory function $D(x)$ is discontinued at the integer arguments $x \ge 1$, the explicit formula is defined by the arithmetic average $D_0(x) = (D(x-\delta) + D(x+\delta))/2$ at the integers values $x \in \mathbb{N}$. Here $\delta > 0$ is an arbitrarily small number.

## 1.2 The Error Terms
The error term is defined by the difference

$$\Delta(x) = \sum_{n \le x} d(n) - \left( \log x + 2\gamma - 1 \right) x .$$

(9)

The generalization to $\Delta_k(x) = \sum_{n \le x} d_k(n) - (P_k(\log x) + 2\gamma - 1)x$, where the residue at $s = 1$ is a polynomial $P_k(x) = \operatorname*{Res}_{s \to 1}\left(\zeta(s)^k x^s / s\right)$ of degree $\deg P_k(x) = k - 1$, and $k \ge 2$ is an integer, is discussed in [IV, p. ]. For example,

$$P_2(x) = x + 2\gamma - 1, \qquad P_3(x) = x^2/2 + (3\gamma - 1)x + 3\gamma^2 - 3\gamma + 3\gamma_1 + 1. \tag{10}$$

The first polynomial $P_2(x)$ is computed in Theorem 3.1, the calculations of the other are similar.

The Dirichlet divisor problem states that $\Delta(x) = O(x^{1/4 + \varepsilon})$ for any arbitrarily small number $\varepsilon > 0$. The works of dozens of authors have reduced the classical estimate $\Delta(x) = O(x^{1/2})$ to approximately $\Delta(x) = O(x^{31})$, see [HL].

Here it is shown that an application of the explicit formula easily produces the standard error term $\Delta(x) = O(x^{1/2})$, and other estimates.

**_Theorem_ 1.2.** For a large number $x \in \mathbb{R} - \mathbb{N}$, the average orders of the divisors counting function has the followings asymptotic formulas.

(i) $\displaystyle \sum_{n \le x} d(n) = -\frac{\pi^2}{12} + (\log x + 2\gamma - 1)x + O(x^{1/2})$,      unconditionally.    (11)

(ii) $\displaystyle \sum_{n \le x} d(n) = -\frac{\pi^2}{12} + (\log x + 2\gamma - 1)x + O(x^{1/4 + \varepsilon})$,     conditional on the RH.

**_Proof_:** (i) Let $\rho = \sigma + it$ denote the zeroes of the zeta function, $1/2 \le \sigma \le 1$, and $t \in \mathbb{R}$. Then, by Theorem 1.1, the estimate can be derived as follows:

$$\begin{aligned} \Delta(x) &= \left| \frac{\pi^2}{3} \sum_{\rho} \frac{\zeta(\rho/2)^2}{\rho \zeta'(\rho)} x^{\rho/2} - \frac{\pi^2}{6} \sum_{n \ge 0} \frac{\zeta(-2n-1)^2}{(2n+1)\zeta'(-2n-1)} x^{-2n-1} \right| \\ &\le \left| \frac{\pi^2}{3} \sum_{\rho} \frac{\zeta(\rho/2)^2}{\rho \zeta'(\rho)} x^{\rho/2} \right| + \left| -\frac{\pi^2}{6} \sum_{n \ge 0} \frac{\zeta(-2n-1)^2}{(2n+1)\zeta'(-2n-1)} x^{-2n-1} \right|, \\ &= O(x^{1/2}). \end{aligned} \tag{12}$$

This follows from Corollary 3.2.

(ii) Assume that the nontrivial zeros of the zeta function are of the form $\rho = 1/2 + it$, $t \in \mathbb{R}$. Proceed as before to determine the error term $\Delta(x) = O(x^{1/4 + \varepsilon})$ for any arbitrarily small number $\varepsilon > 0$.     ■

According to [IR, p. 5], "the conjecture on the upper bound is one of the most difficult problems in analytic number theory, as it does not appear to follow from the Lindelof hypothesis or the Riemann hypothesis."

## 1.3 Omega Result

A result of Hardy shows that

$$\limsup_{x \to \infty} \Delta(x) / x^{1/4} = \infty, \tag{13}$$

and the complementary result of Ingham or Landau, [IV, p. 383], [TN, p. 38], shows that

$$\liminf_{x \to \infty} \Delta(x) / x^{1/4} = -\infty. \tag{14}$$

The proofs of these results, proved by Hardy, Ingham, Landau and other authors, see [KC, p. 195], use elaborate applications of the truncated Voronoi formula

$$\begin{aligned}
\Delta(x) &= \sum_{n \le x} d(n) - \frac{1}{4} - \left( \log x + 2\gamma - 1 \right) x \\
&= \frac{x^{1/4}}{\pi \sqrt{2}} \sum_{n \le N} \frac{d(n)}{n^{3/4}} \cos(4\pi \sqrt{nx} - \pi / 4) + O\left( x^{1/2+\varepsilon} / N^{1/2} \right),
\end{aligned} \tag{15}$$

for $2 \le N < x$ confer [IT, p. 2], [IV, p. 86], [TT]. Some authors have computed quantitative versions such as

$$\Delta(x) = \Omega_{\pm} \left( (x \log x)^{1/4} (\log \log x)^{\alpha} (\log \log \log x)^{-\beta} \right) \tag{16}$$

for some constants $\alpha > 0$, $\beta > 0$, in [SD] and similar references.

The formal explicit formula derived in Theorem 1.1 above provides a different and unified technique for proving various estimates of the summatory function $D(x) = \sum_{n \le x} d(n)$, such as quantitative versions of the limit infimum and limit supremum.

***Theorem* 1.3.** For a large number $x \in \mathbb{R} - \mathbb{N}$, the average order of the divisors counting function has the followings omega asymptotic formula

$$\sum_{n \le x} d(n) = -\frac{\pi^2}{12} + (\log x + 2\gamma - 1)x + \Omega_{\pm}(x^{1/4}). \tag{17}$$

## 1.4 Equivalence of the Divisor Problem and the Circle Problem

The Gauss circle problem states that

$$\sum_{n \leq x} r(n) = \pi x + O(x^{1/4+\varepsilon}),$$ (18)

with $\varepsilon > 0$ arbitrarily small. The Gauss circle problem and the Dirichlet divisor problem are equivalent problem. There is a conversion procedure due to Richert that expresses the circle problem in terms of the divisor problem, refer to [BZ], [CT], [IM], [IV, p. 374] et alii. The Sierpinski formula

$$\sum_{n \leq x} r(n) = \pi x + x^{1/2} \sum_{n \geq 1} \frac{r(n)}{n^{1/2}} J_1(2\pi(nx)^{1/2}),$$ (19)

is similar to the Voronoi formula in (3). But, the summation has a different Bessel function

$$J_v(z) = \sum_{n \geq 0} \frac{(-1)^n}{n! \, \Gamma(v+n+1)} \left( \frac{z}{2} \right)^{v+2n},$$ (20)

which converges for all $0 < |z| < \infty$, and $v \in \mathbb{C}$, see [BZ] for finer details.

An application of the Sierpinski formula of the divisor function reduces the error term to $O(x^{1/3})$. And an unconditional proof of the Gauss circle problem was proposed a few years ago in [CS].

## 2 Generating Functions

The generating function $g(s) = \sum_{n \geq 1} \alpha(n) n^{-s}$ of an arithmetic function $\alpha : \mathbb{N} \rightarrow \mathbb{C}$ determines the shape of the explicit formula

$$\sum_{n \leq x} \alpha(n) = \frac{1}{i2\pi} \int_{c-i\infty}^{c+i\infty} g(s) \frac{x^s ds}{s} = \sum_{s \in Z} \operatorname{Re} s(g,s) + E_\alpha(x),$$ (21)

where $Z = \{ \rho : f(\rho) = \infty \}$ is the set of poles of the function $f(s) = g(s)x^s / s$, and $E_\alpha(x)$ is an error term associated with the arithmetic function $\alpha(n)$, and the integral.

The symbol $\omega(n) = \#\{ p \mid n : p \text{ prime} \}$ denotes the prime divisors counting function, and $\Omega(n) = \#\{ p^\alpha \mid n : p^\alpha \text{ prime power} \}$. The squarefree divisors counting function is defined by $\#\{ d \mid n : \mu(n) \neq 0 \} = 2^{\omega(n)}$.

### 2.1 Some Generating Functions

***Lemma* 2.1.** Let $1 \leq a < q$, and $\gcd(a,q) = 1$. Then, the generating function of the integers $2^{\omega(n)}$ is given by

(i) $\quad \sum_{n \geq 1} \dfrac{2^{\omega(n)}}{n^s} = \dfrac{\zeta(s)^2}{\zeta(2s)},$ \hfill (22)

(ii) $\quad \displaystyle\sum_{n \geq 1, n \equiv a \bmod q} \dfrac{2^{\omega(n)}}{n^s} = \dfrac{\zeta(s)^2}{\zeta(2s)} \prod_{p \mid q} \dfrac{1 - p^{-s}}{1 + p^{-s}}.$

***Proof* :** The function $\omega(p^\alpha q^\beta) = \omega(p^\alpha)\omega(q^\beta)$, $\gcd(p,q) = 1$, is multiplicative but not completely multiplicative, so

$$
\begin{aligned}
\sum_{n \geq 1} \frac{2^{\omega(n)}}{n^s} &= \prod_{p \geq 2}\left(1 + \frac{\omega(p)}{p^s} + \frac{\omega(p^2)}{p^{2s}} + \cdots\right) \\
&= \prod_{p \geq 2}\left(1 + \frac{2}{p^s} + \frac{2}{p^{2s}} + \cdots\right) \\
&= \prod_{p \geq 2}\left(1 + \frac{2}{p^s}\left(1 - \frac{1}{p^s}\right)^{-1}\right) \\
&= \prod_{p \geq 2}\left(1 - \frac{1}{p^s}\right)^{-1}\left(\frac{2}{p^s} + \left(1 - \frac{1}{p^s}\right)\right) \\
&= \prod_{p \geq 2}\left(1 - \frac{1}{p^s}\right)^{-1}\left(1 - \frac{1}{p^s}\right).
\end{aligned}
$$

\hfill (23)

This proves statement (i), to unearth the proof of the second statement (ii), observe that $\sum_{n \equiv a \bmod q} n^{-s} = \zeta(s)\prod_{p \mid q}(1 - p^{-s})$. ∎

Let $n = p_1^{v_1} p_2^{v_2} \cdots p_t^{v_t}$, be the prime decomposition of an integer. For $k \geq 2$, the $k$th divisors function is defined by

$$
d_k(n) = \#\left\{(d_1, \ldots, d_k) : n = d_1 \cdots d_k\right\} = \binom{v_1 + k - 1}{k - 1} \cdots \binom{v_t + k - 1}{k - 1}.
$$

\hfill (24)

This arithmetic function arises as the $n$th coefficient of the $k$ power of the zeta function $\zeta(s)^k$, see [IK, p. 13], and similar references.

***Lemma* 2.2.**     Let $k \geq 2$, let $1 \leq a < q$, and $\gcd(a, q) = 1$. Then, the generating function is given by

(i) $\displaystyle \sum_{n \geq 1} \frac{d_k(n)}{n^s} = \zeta(s)^k$,                                                                   (25)

(ii) $\displaystyle \sum_{n \geq 1, \, n \equiv a \bmod q} \frac{d_k(n)}{n^s} = \zeta(s)^k \prod_{p \mid q} \left(1 - p^{-s}\right)^k$.

There is an extensive literature on these functions, [TT], [IV, p. 352], and other sources.

***Lemma* 2.3.**     The quadratic generating functions are given by

(i) $\displaystyle \sum_{n \geq 1} \frac{d(n^2)}{n^s} = \frac{\zeta(s)^3}{\zeta(2s)}$,                                              (26)

(ii) $\displaystyle \sum_{n \geq 1} \frac{d(n)^2}{n^s} = \frac{\zeta(s)^4}{\zeta(2s)}$.

***Lemma* 2.4.**     The quadratic generating functions are given by

$$\sum_{n \geq 1} \frac{\mu(n) d(n)}{n^s} = \frac{g(s)}{\zeta(s)^2},$$                                                              (27)

where $g(s)$ is analytic for $\Re(s) > 1/2$.

***Proof* :** Expand the infinite sum into a product, then reclassify the different factors:

$$\sum_{n \geq 1} \frac{\mu(n) d(n)}{n^s} = \prod_{p \geq 2} \left(1 - \frac{2}{p^s}\right)$$

$$= \prod_{p \geq 2} \left(1 - \frac{1}{p^s}\right)^2 \left(1 - \frac{1}{p^{2s}}\left(1 - \frac{1}{p^s}\right)^{-2}\right)$$                (28)

$$= \frac{g(s)}{\zeta(s)^2}.$$

Since the second term is analytic for $\Re(s) > 1/2$, it proves statement.     ∎

The divisors function $d(n) = \sum_{d \mid n} 1$ is a special case of the more general divisors functions $\sigma_a(n) = \sum_{d \mid n} d^a$, and its generating function $g(s) = \sum_{n \geq 1} \sigma_a(n) n^{-s} = \zeta(s)\zeta(s - a)$. The generating functions of various arithmetic functions are linked via some parameters. For example, the proof of (ii) can be extracted from the more general generating function of product of divisor functions

$$\sum_{n \geq 1} \frac{\sigma_a(n)\sigma_b(n)}{n^s} = \frac{\zeta(s)\zeta(s-a)\zeta(s-b)\zeta(s-a-b)}{\zeta(2s-a-b)},$$  (29)

for $a, b \in \mathbb{C}$. As another example, at $k = -1$, the expression in Lemma 2.1-i yields the generating function $\sum_{n \geq 1} \mu(n)n^{-s} = \zeta(s)^{-1}$ of the Mobius function $\mu(n) \to \{-1, 0, 1\}$.

Some generalization of generating functions to quadratic fields are studied in [KN, p. 2], [GN], et cetera.

### 2.2 Level $k$ Reduction Identities

The level $k \geq 2$ reduction identities express various complicated divisors function in terms of simpler arithmetic functions.

***Lemma* 2.5.**    Let $k \geq 2$, let $d_k(n) = \#\{(d_1,...,d_k) : d_1 \cdots d_k = n \}$, and $d(n) = d_2(n)$. Then,

(i)   $2^{\omega(n)} = \sum_{d \mid n} \mu^2(d)$,  (30)

(ii)   $d(n^2) = \sum_{d^2 \mid n} \mu(d)d_3(n/d^2)$,

(iii)   $d(n)^2 = \sum_{d^2 \mid n} \mu(d)d_4(n/d^2)$.

The first formula (i) tallies the number of squarefree divisors of $n \geq 1$.

***Lemma* 2.6.**    Let $k \geq 2$, let $d_k(n) = \#\{(d_1,...,d_k) : d_1 \cdots d_k = n \}$, and $d(n) = d_2(n)$. Then it has the recursive formula

$$d_k(n) = \sum_{d \mid n} d_{k-1}(n).$$  (31)

***Proof*:** Compute the $n$th coefficient of the $k$ power of the zeta function $\zeta(s)^k$ recursively:

$$\sum_{n \geq 1} \frac{d_k(n)}{n^s} = \zeta(s) \cdot \zeta(s)^{k-1} = \sum_{n \geq 1} \frac{1}{n^s} \cdot \sum_{n \geq 1} \frac{d_{k-1}(n)}{n^s}.$$  (32)

**Problems.**

1. Let $q > 1$ be a modulo of an arithmetic progression. Use Lemma 2.1-ii to compute the explicit formula for

$$\sum_{n \leq x, \, n \equiv a \bmod q} 2^{\omega(n)} = \frac{f(x)}{\varphi(q)} + \sum_{\rho} g(\rho) x^{\rho} + \sum_{n < 0} h(n) x^{-n}, \tag{33}$$

where $f$, $g$, $h$ are some continuous functions, and the index $\rho$ runs over the nontrivial zeroes of the zeta function $\zeta(s)$.

2. Use the sigma-phi identity $\left(\sigma(n)/n\right)\left(\varphi(n)/n\right) = \prod_{p^{\varepsilon} \| n}\left(1 - p^{-\alpha - 1}\right)$ or some other identity to modify the generating function $g(s) = \sum_{n \geq 1} \sigma(n) n^{-s} = \zeta(s)\zeta(s - 1)$, and use it to construct an explicit formula for

$$\sum_{n \leq x} \sigma(n) = f(x) + \sum_{\rho} g(\rho) x^{\rho} + \sum_{n < 0} h(n) x^{-n}, \tag{34}$$

where $f$, $g$, $h$ are some continuous functions, and the index $\rho$ runs over the nontrivial zeroes of the zeta function $\zeta(s)$.

# 3 Some Explicit Formulas

The domain of definition of the generating function $\sum_{n \geq 1} 2^{\omega(n)} n^{-s} = \zeta(s)^2 / \zeta(2s)$ of squarefree integers has a visible dependence on the zeroes and pole of the zeta function. The theory of residue in tandem with the Perron formula transform this information into a nice explicit formula for the partial sum $\sum_{n \leq x} 2^{\omega(n)}$.

***Theorem* 3.1.** For a large number $x \in \mathbb{R} - \mathbb{N}$, the explicit formula of the average order $\sum_{n \leq x} 2^{\omega(n)}$ of the squarefree divisors counting function $2^{\omega(n)}$ is given by

$$\sum_{n \leq x} 2^{\omega(n)} = -\frac{1}{2} + \frac{6}{\pi^2}\left(\log x + 2\gamma - 1 - 12\zeta'(2)/\pi^2\right)x + 2\sum_{\rho} \frac{\zeta(\rho/2)^2}{\rho\zeta'(\rho)} x^{\rho/2}$$
$$- \sum_{n \geq 0} \frac{\zeta(-2n-1)^2}{(2n+1)\zeta'(-2n-1)} x^{-2n-1}, \tag{35}$$

where $\rho$ ranges over the nontrivial zeroes of the zeta function $\zeta(s)$.

***Proof***: Let $Z = \left\{\, 0, \, 1, \, -n, \, \rho/2 : \zeta(\rho) = 0, \text{ and } n \geq 1 \,\right\}$ be the set of poles of the function $f(s) = \zeta(s)^2 x^s / (s\zeta(2s))$, and evaluate the Perron integral:

$$\sum_{n \le x} 2^{\omega(n)} = \int_{c-i\infty}^{c+i\infty} \frac{\zeta(s)^2}{\zeta(2s)} \frac{x^s ds}{s} = \sum_{s \in Z} \operatorname{Re} s(f,s), \tag{36}$$

where $c > 1$ is a constant, see [IK, p. 151]. The contributions by the top, bottom, and left line integrals are negligible, and are not shown. The contributions of these integrals are small because of the exponential decay of the function $f(s)$ on the regions $\Re e(s) < 0$ and $\Re e(s) > 1/2$. The sum of residues on the right side of the line integral (36) involves four types of residues:

**(i).** The residue of $f(s) = \zeta(s)^2 x^s /(s\zeta(2s))$ at the simple pole at $s = 0$.

$$\operatorname{Re} s(f, s = 0) = \lim_{s \to 0} s \cdot \frac{\zeta(s)^2 x^s}{s\zeta(2s)} = \lim_{s \to 0} \frac{\zeta(s)^2 x^s}{\zeta(2s)} = \frac{(-1/2)^2}{-1/2} = -\frac{1}{2}. \tag{37}$$

**(ii).** The residue of $f(s) = \zeta(s)^2 x^s /(s\zeta(2s))$ ) at the simple poles at $s = \rho/2$, such that $\zeta(\rho) = 0$.

$$\operatorname{Re} s(f, s = \rho/2) = \lim_{s \to \rho/2} (s - \rho/2) \cdot \frac{\zeta(s)^2 x^s}{s\zeta(2s)} = \lim_{s \to \rho/2} \frac{\zeta(s)^2 x^s}{2s\zeta'(2s)} = \frac{\zeta(\rho/2)^2 x^{\rho/2}}{\rho\zeta'(\rho)}. \tag{38}$$

This is done via the l'Hospital rule since the Laurent series of the zeta function at $s = \rho/2$ is unknown.

**(iii).** The residue of $f(s) = \zeta(s)^2 x^s /(s\zeta(2s))$ ) at the simple pole at $s = 1$.

$$\begin{aligned}
\operatorname{Re} s(f, s_1) &= \frac{d^{(m-1)}}{ds}(s - s_1)^m f(s) \bigg|_{s=s_1} = \frac{d}{ds}(s-1)^2 \frac{\zeta(s)^2 x^s}{s\zeta(2s)} \bigg|_{s=1} \\
&= \frac{6}{\pi^2} \big( \log x + (2\gamma - 1 - 12\zeta'(2)/\pi^2) \big) x.
\end{aligned} \tag{39}$$

This is computed using the expansion

$$\begin{aligned}
(s-1)^2 f(s) &= \frac{(s-1)^2 x^s}{s\zeta(2s)} \zeta(s)^2 \\
&= \frac{(s-1)^2 x^s}{s\zeta(2s)} \left( \frac{1}{s-1} + \gamma + \gamma_1(s-1) + \gamma_2(s-1)^2 + \cdots \right)^2 \\
&= \frac{(s-1)^2 x^s}{s\zeta(2s)} \left( \frac{1}{(s-1)^2} + \frac{2\gamma}{s-1} + \gamma^2 + 2\gamma_1 + \cdots \right) \\
&= \frac{x^s}{s\zeta(2s)} \big( 1 + 2\gamma(s-1) + (\gamma^2 + 2\gamma_1)(s-1)^2 + \cdots \big),
\end{aligned} \tag{40}$$

see [DL, eq. 25.2.4] for more information on the Laurent series of the zeta function at $s = 1$.

**(iv).** The residue of $f(s) = \zeta(s)^2 x^s / (s\zeta(2s))$ at the simple poles at $s = -n$, where $n \geq 1$.

$$\operatorname{Re} s(f, s = -n) = \lim_{s \to -n} (s + n) \cdot \frac{\zeta(s)^2 x^s}{s\zeta(2s)} = \lim_{s \to -n} \frac{\zeta(s)^2 x^s}{2s\zeta'(2s)} = \frac{\zeta(-n)^2 x^{-n}}{-2n\zeta'(-2n)}. \tag{41}$$

This is done via the l'Hospital rule since the Laurent series of the zeta function at $s = -n$ is unknown.

Replacing these data into the Perron formula, and using $\zeta(-n) = 0$ at the even integers $-n = -2m$, $m \geq 1$, complete the evaluation. ∎

**Corollary 3.2.** For a large number $x \in \mathbb{R} - \mathbb{N}$, the explicit formula of the average order $\sum_{n \leq x} 2^{\omega(n)} / n$ of the squarefree divisors counting function $2^{\omega(n)}$ is given by

$$\sum_{n \leq x} \frac{2^{\omega(n)}}{n} = \frac{6}{\pi^2} \left( \frac{\log^2 x}{2} + (2\gamma - 12\zeta'(2)/\pi^2) \log x \right) + 2\gamma - 1 + 2\sum_{\rho} \frac{\zeta(\rho/2)^2}{\rho\zeta'(\rho)} x^{\rho/2-1}$$
$$- \sum_{n \geq 0} \frac{\zeta(-2n-1)^2}{(2n+1)\zeta'(-2n-1)} x^{-2n-2}, \tag{42}$$

where $\rho$ ranges over the nontrivial zeroes of the zeta function $\zeta(s)$.

**Corollary 3.3.** Let $\rho = \sigma + it$ be the nontrivial the zeroes of the zeta function $\zeta(s)$. Then,

$$\left| \sum_{\rho} \frac{\zeta(\rho/2)^2}{\rho\zeta'(\rho)} \right| \leq O(1). \tag{43}$$

**Proof**: For a large number $x \in \mathbb{R} - \mathbb{N}$, the error term of the explicit formula in Theorem 3.1 is the expression

$$\Delta_0(x) = \sum_{n \leq x} 2^{\omega(n)} - \frac{6}{\pi^2} \left( \log x + 2\gamma - 1 - 12\zeta'(2)/\pi^2 \right) x$$
$$= 2\sum_{\rho} \frac{\zeta(\rho/2)^2}{\rho\zeta'(\rho)} x^{\rho/2} - \sum_{n \geq 0} \frac{\zeta(-2n-1)^2}{(2n+1)\zeta'(-2n-1)} x^{-2n-1}. \tag{44}$$

Here, the absolute value of the second infinite sum on the right is bounded by a constant, and the absolute value of the first infinite sum dominates the expression on the right side. Hence, it satisfies

$$\left| \Delta_0(x) \right| \le C \left| 2\sum_{\rho} \frac{\zeta(\rho/2)^2}{\rho\zeta'(\rho)} x^{\rho/2} \right| = O(x^{1/2}),$$

(45)

where $C > 0$ is a constant. Next, let $\rho = \sigma + it$, where $1/2 \le \sigma \le 1$, and $t \in \mathbb{R}$. This proves the claim. ∎

**Remark 1.** A different analysis of the infinite component, using separate estimates of the individuals functions $\zeta(s)$, $\zeta'(s)$ on the critical strip $0 < \Re e(s) < 1$ is of independent interest, but could be difficult.

## 4 Summatory Functions of Squarefree and Nonsquarefree Divisors

Various estimates of the summatory functions $\sum_{n \le x} 2^{\omega(n)}$ and $\sum_{n \le x} 2^{\Omega(n)}$ of $2^{\omega(n)}$ and $2^{\Omega(n)}$ respectively are presented in this Section. The counting functions $2^{\omega(n)} = \#\{ d \mid n : \mu(n) \ne 0 \} = 2^{\Omega(n)}$ coincide if $n$ is a squarefree integer. Otherwise, for a nonsquarefree integer $n$, there is an strict inequality $2^{\omega(n)} < d(n) < 2^{\Omega(n)}$.

***Lemma* 4.1.** For a large number $x \in \mathbb{R}$, the average order of $2^{\omega(n)}$ is given by

(i) $\sum_{n \le x} 2^{\omega(n)} = \frac{6}{\pi^2} \big( \log x + 2\gamma - 1 - 2\zeta'(2)/\zeta(2) \big) x + O(x^{1/2} e^{-c\sqrt{\log x}})$.

(46)

(ii) $\sum_{n \le x} \frac{2^{\omega(n)}}{n} = \frac{6}{\pi^2} \left( \frac{\log^2 x}{2} + (2\gamma - 2\zeta'(2)/\zeta(2)^2)\log x \right) + 2\gamma - 1 + O(x^{-1/2} e^{-c\sqrt{\log x}})$.

***Proof* :** (i) Consider the identity

$$\sum_{n \le x} 2^{\omega(n)} = \sum_{n \le x} \sum_{d^2 m = n} \mu(d) d(m) = \sum_{n \le x} \sum_{d^2 \mid n, n \le x} \mu(d) d(n/d^2).$$

(47)

Reversing the order of summation and evaluating the inner finite sum returns

$$\sum_{n \leq x} 2^{\omega(n)} = \sum_{d \leq x} \mu(d) \sum_{n \leq x/d^2} d(n)$$

$$= \sum_{d \leq x} \mu(d) \left( \frac{x}{d^2} (\log(x/d^2) + 2\gamma - 1) + O(\frac{x^{1/2}}{d}) \right) \qquad (48)$$

$$= x \log x \sum_{d \leq x} \frac{\mu(d)}{d^2} - x \sum_{d \leq x} \frac{\mu(d) \log d^2}{d^2} + (2\gamma - 1) x \sum_{d \leq x} \frac{\mu(d)}{d^2} + O(x^{1/2} \sum_{d \leq x} \frac{\mu(d)}{d}) \,.$$

The second term in the previous expression can be transformed to

$$x \sum_{d \leq x} \frac{\mu(d) \log d^2}{d^2} = 2x \sum_{d \geq 1} \frac{\mu(d) \log d}{d^2} - 2x \sum_{d > x} \frac{\mu(d) \log d}{d^2} = 2x \frac{\zeta'(2)}{\zeta(2)^2} + O(\log x) \,, \qquad (49)$$

and

$$\sum_{d \leq x} \frac{\mu(d)}{d} = O(e^{-c\sqrt{\log x}}) \,, \qquad (50)$$

for some absolute constant $c > 0$. Returning these into the previous expression, and simplifying the other terms yield (i).

(ii) The second finite sum can be derived by partial summation. ∎

A much earlier analysis of the above result, but more complicated appears in [GD], see also [MV, p. 42]. It should be observe that the calculation of the upper bound (50) is not very easy, however, the weaker upper bound

$$\left| \sum_{d \leq x} \frac{\mu(d)}{d} \right| < 1 \,, \qquad (51)$$

is quite easy.

***Lemma* 4.2.**     Assume the RH. Then, for a large number $x \in \mathbb{R}$, the average order of $2^{\omega(n)}$ is given by

(i)   $\sum_{n \leq x} 2^{\omega(n)} = \frac{6}{\pi^2} \left( \log x + 2\gamma - 1 - 2\zeta'(2)/\zeta(2) \right) x + O(x^{1/4+\varepsilon}) \,.$ \qquad (52)

(ii)  $\sum_{n \leq x} \frac{2^{\omega(n)}}{n} = \frac{6}{\pi^2} \left( \frac{\log^2 x}{2} + (2\gamma - 1 - 2\zeta'(2)/\zeta(2)^2) \log x \right) + 2\gamma - 1 + O(x^{-3/4+\varepsilon}) \,.$

***Proof* :** (i) There are two proofs. One uses the same analysis as (35) and the second uses the explicit formula in Theorem 3.1. (ii) The second finite sum can be derived from (i) by partial summation. ∎

***Lemma* 4.3.** For a large number $x \in \mathbb{R}$, the average order of $2^{\Omega(n)}$ is given by

(i) $\displaystyle \sum_{n \leq x} 2^{\Omega(n)} = a_0 x \log^2 x + O(x \log x)$, $\hspace{3cm}$ (53)

(ii) $\displaystyle \sum_{n \leq x} \frac{2^{\Omega(n)}}{n} = a_0 \log^2 x + O(\log x)$,

where $a_0 = (8 \log 2)^{-1} \prod_{p > 2} (1 + p^{-1}(p-2)^{-1})$ is a constant.

***Corollary* 4.4.** For a large number $x \in \mathbb{R}$, the average order of the nonsquarefree divisors is given by

$$\sum_{n \leq x, \, \mu(n)=0} d(n) = -\frac{2\zeta'(2)}{\zeta(2)^2} x + O(x^{1/2} e^{-c\sqrt{\log x}}),$$ $\hspace{2cm}$ (54)

for some absolute constant $c > 0$.

# 5 Omega Result

***Theorem* 5.1.** For a large number $x \geq 1$, the average orders of the squarefree divisor function satisfies the omega relation

$$\sum_{n \leq x} 2^{\omega(n)} = -\frac{1}{2} + \frac{6}{\pi^2} \left( \log x + 2\gamma - 1 - 12\zeta'(2)/\pi^2 \right) x + \Omega_{\pm}(x^{1/4}).$$ $\hspace{1cm}$ (55)

***Proof***: For each nontrivial zero $\rho = 1/2 + it$, $t \in \mathbb{R}$, the generating function $\sum_{n \geq 1} 2^{\omega(n)} n^{-s} = \zeta(s)^2/\zeta(2s)$ has a singularity at $\rho/2$. Therefore, the difference

$$\sum_{n \leq x} 2^{\omega(n)} - \frac{6}{\pi^2} \left( \log x + 2\gamma - 1 - 12\zeta'(2)/\pi^2 \right) x = -\frac{1}{2} + 2 \sum_{\rho} \frac{\zeta(\rho/2)^2}{\rho \zeta'(\rho)} x^{\rho/2}$$

$$- \sum_{n \geq 0} \frac{\zeta(-2n-1)^2}{(2n+1)\zeta'(-2n-1)} x^{-2n-1},$$ $\hspace{1cm}$ (56)

where $\rho$ ranges over the nontrivial zeros of the zeta function $\zeta(s)$, assumes arbitrarily large values at any nontrivial zero $\rho/2 = 1/4 + it/2$. This in turns implies that the absolute value of the difference satisfies

$$\left| \frac{1}{2} + 2\sum_{\rho} \frac{\zeta(\rho/2)^2}{\rho\zeta'(\rho)} x^{\rho/2} + \sum_{n \geq 0} \frac{\zeta(-2n-1)^2}{(2n+1)\zeta'(-2n-1)} x^{-2n-1} \right| \geq x^{1/4} \left| \sum_{\rho} \frac{\zeta(\rho/2)^2}{\rho\zeta'(\rho)} + O\left(\frac{1}{x^{1/4}}\right) \right|^{1/4} \tag{57}$$
$$\geq c_0 x^{1/4},$$

where $c_0 > 0$ is a constant. ∎

# 6. The Divisor Function

The divisor function $d(n) = \sum_{d \mid n} 1$ tallies the number of divisors of a natural number $n \in \mathbb{N}$. Given the prime decomposition of the integer $n = p_1^{v_1} \cdot p_2^{v_2} \cdots p_t^{v_t}$, where the $p_i$ are primes, and $v_i \geq 1$, it can be computed using the formula

$$d(n) = \prod_{1 \leq i \leq t} \left(1 + v_i\right). \tag{58}$$

The value $d(n)$ always falls within the range $2^{\omega(n)} \leq d(n) \leq 2^{\Omega(n)}$, where $\omega(n) = \#\{p \mid n : p \text{ is prime}\} = t$, and $\Omega(n) = v_1 + \cdots + v_t$ counts the total number of prime powers divisors, including multiplicities.

The divisor function is a quasirandomly varying function of the integer $n \geq 1$. It can unpredictably jump from its minimal value $d(n) = 2$ at the primes arguments $n = p_1^{v_1} \cdot p_2^{v_2} \cdot p_3^{v_3} \cdots p_t^{v_t} - 1$ to its maximal value

$$d(n) = (v_1 + 1)(v_2 + 1)\cdots(v_t + 1) \geq 2^t \tag{59}$$

at the highly composite integers $n = p_1^{v_1} \cdot p_2^{v_2} \cdots p_t^{v_t}$, with $v_1 \geq v_2 \geq \cdots \geq v_t \geq 1$, $t \geq 1$, and vice versa. The upper bound is often given in the form $d(n) = O(n^\varepsilon)$, with $\varepsilon > 0$ arbitrarily small.

*Lemma* 6.1.  (Hermite formula) Let $n \in \mathbb{N}$ be an integer. Then,

$$d(n) = 2 \sum_{d \mid n, n \leq \sqrt{x}} 1 - \delta(n), \tag{60}$$

where the delta function is defined by

$$\delta(n) = \begin{cases} 1 & \text{if } n = m^2, \\ 0 & \text{if } n \neq m^2. \end{cases} \tag{61}$$

***Example*** **1.**  For the integers $n = 100$, and $n = 101$, the number of divisors are respectively:

$$d(100) = 2 \sum_{d \mid n, \, n \leq \sqrt{100}} 1 \; - \delta(100) = 9, \quad \text{and} \quad d(101) = 2 \sum_{d \mid n, \, n \leq \sqrt{101}} 1 \; - \delta(101) = 2. \tag{62}$$

The limit supremum has the following form, [HW, p. 345], [NW, p. 254]:

***Theorem*** **6.2.**  (Wigert)   Let $n \in \mathbb{N}$ be an integer, and let $d(n) = \sum_{d \mid n} 1$ be the divisor function. Then

$$\limsup_{n \to \infty} \frac{\log d(n) \log \log n}{\log n} = \log 2. \tag{63}$$

***Theorem*** **6.3.**   Let $\varepsilon > 0$ be a small real number. Then
(i)   $d(n) < 2^{(1+\varepsilon) \log n / \log \log n}$,                              for all integers $n \geq n_\varepsilon$.
(ii)   $d(n) > 2^{(1-\varepsilon) \log n / \log \log n}$,                              for infinitely many integers $n \geq n_\varepsilon$.

These results are standard in the literature, the proofs appear in [NV, p. 398] and similar sources.

***Theorem*** **6.4.**   ([HW, p. 478])   Let $\varepsilon > 0$ be an arbitrarily small number. Then, the inequality

$$2^{(1-\varepsilon) \log n / \log \log n} < d(n) < 2^{(1+\varepsilon) \log n / \log \log n} \tag{64}$$

holds for almost all integers $n \in \mathbb{N}$.

***Theorem*** **6.5.**      The ratio $d(n) / 2^{\log n / \log \log n}$ is dense in $[\, 0, \infty \,]$.

***Proof***: Let $\varepsilon > 0$ be an arbitrarily small number. By Theorem 6.4, there exists infinitely many integers $n \in \mathbb{N}$ such that

$$2^{-\varepsilon \log n / \log \log n} < \frac{d(n)}{2^{\log n / \log \log n}} < 2^{\varepsilon \log n / \log \log n}. \tag{65}$$

Since the number $\varepsilon > 0$ is arbitrary, there are arbitrary values of the ratio in the interval $[\, 2^{-\varepsilon \log n / \log \log n}, 2^{\varepsilon \log n / \log \log n} \,] \subset [\, 0, \infty \,]$.                              ∎

## 7 Simple Estimates of the Average Order

The average order $\sum_{n \leq x} d(n)$ is a much better estimator of the behavior of the values of the divisor function $d(n) = \sum_{d \mid n} 1$ as the integer $n$ varies over the interval $[\ 1,\ x\ ]$. A few computations of the average order of the divisors function are sketched in this Section.

The first calculation gives the precise main term, and rough estimate of the error term.

***Lemma* 7.1.** If $x \geq 1$ is a sufficiently large real number, then

$$\sum_{n \leq x} d(n) = x \log x + O(x) \ . \tag{66}$$

***Proof***: Use an inversion of the order of summation to write the finite sum as follows:

$$\sum_{n \leq x} d(n) = \sum_{n \leq x,\ d \mid n} 1 = \sum_{d \leq x,\ n \leq x, d \mid n} 1 = \sum_{d \leq x} \left[ \frac{x}{d} \right], \tag{67}$$

where $[\ x\ ]$ largest integer function. Substitute the integer/fractional part identity $[\ x\ ] = x - \{\ x\ \}$ to express the inner finite sum as

$$\begin{aligned}
\sum_{n \leq x} d(n) &= \sum_{d \leq x} \left( \frac{x}{d} - \{\ x / d\ \} \right) \\
&= x \sum_{d \leq x} \frac{1}{d} - \sum_{d \leq x} \{\ x / d\ \}.
\end{aligned} \tag{68}$$

Next, use a rough estimate of the harmonic finite sum $\sum_{n \leq x} 1/n = \log x + O(1)$ to complete the proof. ∎

### 7.1 A Sharp Estimate of the Average Order

The second calculation gives the precise main term, and a significantly sharper estimate of the error term.

***Lemma* 7.2.** If $x \geq 1$ is a sufficiently large real number, then

$$\sum_{n \leq x} d(n) = \left( \log x + 2\gamma - 1 \right) x + O(x^{1/2}) \ . \tag{69}$$

**Proof**: Use an inversion of the order of summation to write the finite sum as follows:

$$\sum_{n \leq x} d(n) = \sum_{n \leq x,} \sum_{d \mid n} 1 = \sum_{d \leq x,} \sum_{n \leq x, d \mid n} 1 = \sum_{d \leq x} \left[ \frac{x}{d} \right], \tag{70}$$

where $[x]$ is largest integer function. Use the integer/fractional part identity $[\, x \,] = x - \{\, x \,\}$ to express the inner finite sum as

$$\sum_{n \leq x} d(n) = \sum_{d \leq x} \left( \frac{x}{d} - \{\, x/d \,\} \right)$$
$$= x \sum_{d \leq x} \frac{1}{d} - \sum_{d \leq x} \{\, x/d \,\}. \tag{71}$$

A precise estimate of the harmonic finite sum has the asymptotic formula

$$\sum_{n \leq x} \frac{1}{n} = \log x + \gamma + O(1/x), \tag{72}$$

see [VL], and the finite sum of fractional parts has the asymptotic formula

$$\sum_{n \leq x} \{\, x/n \,\} = (1 - \gamma)x + O(x^{1/2}) \tag{73}$$

see [PL, p. ]. Substituting these data lead to

$$\sum_{n \leq x} d(n) = x \left( \log x + \gamma + O(1/x) \right) - \left( (1-\gamma)x + O(x^{1/2}) \right)$$
$$= x \log x + (2\gamma - 1)x + O(x^{1/2}), \tag{74}$$

where $\gamma = .5772...$ is the Euler constant. ∎

The standard references for these calculations are [AP, p. 57], [MV, p. 38].

## 7.3 A Sharp Estimate of the Average Order

The third approach to an estimate of the average order of the divisor function uses a lattice points counting approach in the calculations. This method, known as the hyperbola method, gives the precise main term, and a significantly sharper estimate of the error term.

The average order $\sum_{n \leq x} d(n)$ of the divisors counting function $d(n) = \sum_{d \mid n} 1$ can be viewed as the number of integral points under the hyperbola $uv = x$. This done in [HW, p.349], [BR,

Lemma 3.1], [DF, p. 58], [IM, p. 63], [TN, p. 38].

**Lemma 7.3.** If $x \geq 1$ is a sufficiently large real number, then

$$\sum_{n \leq x} d(n) = x \log x + (2\gamma - 1)x + O(x^{1/2}) \ . \tag{75}$$

**Proof**: The number of lattice points up to $x$, and on the hyperbola $uv = n \leq x$ is as follows:

$$\sum_{n \leq x} d(n) = 2 \sum_{n \leq x^{1/2}} \left( [x/n] - n \right) + \left[ x^{1/2} \right]$$

$$= 2 \sum_{n \leq x^{1/2}} \left( \frac{x}{n} - n \right) - 2 \sum_{d \leq x^{1/2}} \{ x/d \} + \left[ x^{1/2} \right] \tag{76}$$

where $[x]$ is the largest integer function, and $\{x\} = x - [x]$ is the fractional part function. Evaluating it returns

$$\sum_{n \leq x} d(n) = 2x \sum_{n \leq x^{1/2}} \frac{1}{n} - 2 \sum_{n \leq x^{1/2}} n - 2 \sum_{n \leq x^{1/2}} \{ x/n \} + \left[ x^{1/2} \right]$$

$$= 2x \left( \log x^{1/2} + \gamma + O(x^{-1/2}) \right) - \left( \left[ x^{1/2} \right] \left( \left[ x^{1/2} \right] + 1 \right) \right) - 2 \sum_{n \leq x^{1/2}} \{ x/n \} + \left[ x^{1/2} \right] \tag{77}$$

$$= x \log x + (2\gamma - 1)x + O(x^{1/2}) ,$$

where $\gamma = .5772...$ is the Euler constant. ∎

The concept of counting points within a bounded region is generalized to other curves, in [HL, p. ] and related references. The lattice points method for counting points inside the quadratic curves $Q(u,v) = au^2 + buv + cv^2$, which generalizes the hyperbola $Q(u,v) = uv$ associated with the divisors function, are carried out in [NR].

**Lemma 7.4.** If $x \geq 1$ is a sufficiently large real number, then

$$\sum_{n \leq x} \frac{d(n)}{n} = \frac{1}{2} \log^2 x + (2\gamma - 1) \log x + a_1 + O(x^{-2/3}) , \tag{78}$$

where $a_1$ is a constant.

Use partial summation, and Lemma 6.3. The proof of this finite sum appears in [BR, Corollary 2.2].

**Lemma 7.5.**   If $x \geq 1$ is a sufficiently large real number, then

$$\sum_{n \leq x} d(n^2) = ax \log^2 x + bx \log x + cx + O(x^{1/2} e^{-c\sqrt{\log x}}),  \tag{79}$$

where $a$, $b$, $c > 0$ are constants.

**Proof**: Use the generating function $\sum_{n \geq 1} d(n^2) n^{-s} = \zeta(s)^3 \zeta(2s)^{-1}$ in Lemma 2.3, other details are given in [IS, p. 4]. ∎

**Lemma 7.6.**   Assume the RH. If $x \geq 1$ is a sufficiently large real number, then

$$\sum_{n \leq x} d(n^2) = ax \log^2 x + bx \log x + cx + O(x^{1/4+\varepsilon}),  \tag{80}$$

where $a$, $b$, $c > 0$ are constants, and $\varepsilon > 0$ is an arbitrarily small number.

**Lemma 7.7.**   (Ramanujan)   If $x \geq 1$ is a sufficiently large real number, then

$$\sum_{n \leq x} d(n)^2 = Ax \log^3 x + Bx \log^2 x + Cx \log x + O(x^{3/5+\varepsilon}),  \tag{81}$$

where $A$, $B$, $C > 0$ are constants.

**Proof**: A few recent proofs are given in [RS], [JS], and other authors. ∎

## 7.4 The Voronoi Estimate of the Average Order
The fourth approach is based on the Voronoi formula, this analysis is far more involved, and the reader is referred to the literature for the proof.

**Lemma 7.8.**   If $x \geq 1$ is a sufficiently large real number, then

$$\sum_{n \leq x} d(n) = x \log x + (2\gamma - 1)x + O(x^{1/3}).  \tag{82}$$

## 8 Binary Divisors Problems

The binary divisors problems deals with the convolutions of the divisors function. For a fixed integer $k \geq 2$, and $m \in \mathbb{Z}$, let

$$\sum_{n \leq x} d_k(n) d_k(n+m) = x P_{2k-2}(\log x, m) + \Delta_k(x, m). \qquad (83)$$

Here $P_{2k-2}(z, m) \in \mathbb{R}[z]$ is a polynomial with real coefficients, of degree $\deg(P_{2k-2}) = 2k - 2 \geq 2$, and $\Delta_k(x, m)$ is the error term.

This case has the polynomial $P_2(z, m) = az^2 + bz + c$, where $\zeta(2)^{-1} \sigma(m)/m$. Beside this, very few results have been achieved on this area of analytic number theory.

***Theorem* 8.1.** (Ingham) Let $m \geq 1$ be fixed, and let $x \geq 1$ is a sufficiently large real number, then

(i) $\quad \displaystyle\sum_{n \leq x} d(n) d(n - m) = \frac{6}{\pi^2} \sigma(m) \log^2 x + O(\log x)$, $\qquad (84)$

(ii) $\quad \displaystyle\sum_{n \leq x} d(n) d(n + m) = \frac{6}{\pi^2} \frac{\sigma(m)}{m} x \log^2 x + O(x \log x)$.

where $a$, $b$, $c > 0$ are constants, and $\varepsilon > 0$ is an arbitrarily small number.

The proofs of these results are essentially the enumerations of the number of integer solutions of the equation $m = ab \pm cd$, $n = c$, where $a$, $b$, $c$, $d \in \mathbb{N}$. However, the analysis is lengthy. The binary divisor problem is discussed in [BB].

## Problems.

1. Let $m \geq 1$ be a fixed integer. Determine an asymptotic formula for

$$\sum_{n \leq x} d_k(n) d_k(n + m) \qquad (85)$$

for some $k \geq 3$.

2. Let $q > 1$ be a modulo of an arithmetic progression, and $m \geq 1$ a fixed integer. Determine an asymptotic formula for

$$\sum_{n \leq x, \, n \equiv a \bmod q} d_k(n) d_k(n + m) \qquad (86)$$

for some $k \geq 2$.

## 9. The Divisor Function Over Arithmetic Progressions

Some of the works on the divisor summatory function over arithmetic progressions was done in [SS], [VC], and [NW]. The divisor function over an arithmetic progression $\{ n \equiv a \bmod q : gdc(a, q) = 1 \}$ is defined by $d(n, q, a) = \#\{ d \mid n : d \equiv a \bmod q \}$.

The formulas for computing this function over arithmetic progressions are more complex than for the plain divisor function $d(n) = \#\{ d \mid n \}$. For example, let $n = p_1^{u_1} p_2^{u_2} \cdots p_k^{u_k} q_1^{v_1} q_2^{v_2} \cdots q_m^{v_m}$, $p_i \equiv 3 \bmod 4$, $q_j \equiv 1 \bmod 4$ be the prime decomposition of an integer, then for the residue class $\{ n \equiv 1 \bmod 4 \}$, the formula has the shape

$$d(n, 4, 1) = \begin{cases} 1 & \text{if } n \neq a^2 + b^2, \\ d(n) & \text{if } n = a^2 + b^2, \text{ and } n = p_1^{u_1} p_2^{u_2} \cdots p_k^{u_k}, \\ d(n) - d_3(n) & \text{if } n = a^2 + b^2, \text{ and } n \neq p_1^{u_1} p_2^{u_2} \cdots p_k^{u_k}, \end{cases} \tag{87}$$

where $d_3(n) = \#\{ d \mid n : d \equiv 3 \bmod 4 \}$ is a correction factor. Similarly, for the residue class $\{ n \equiv 3 \bmod 4 \}$, the formula has the more complicated formula.

The corresponding summatory function is

$$\sum_{n \leq x} d(n, q, a) \ . \tag{88}$$

**Theorem** 9.1. ([NW])  Let $x \in \mathbb{R}$, be a large number, and let $a, q \in \mathbb{N}$ be integers such that $\gcd(a, q) = 1$. Then, there exists a constant $\theta < 1/3$ such that the asymptotic formula

$$\sum_{n \leq x} d(n, q, a) = \frac{1}{q} x \log x + (\gamma(a, q) - \frac{1 - \gamma}{q}) x + O((\frac{x}{q})^\theta) , \tag{89}$$

holds uniformly in $a$, $q$ and $x$ provided that $a < q \leq x$.

The error term has an analogous expression as

$$\Delta(x, q, a) = \sum_{n \leq x} d(n, q, a) - \frac{1}{q} \big( \log x + 2\gamma - 1 \big) x - (\gamma(a, q) - \frac{1 - \gamma}{q}) x = O((\frac{x}{q})^\theta) \tag{90}$$

for some constant $\theta < 1/3$.

Some specific cases are known in the literature. For, example, in [CT] it is shown that

(i) $\displaystyle\sum_{n \le x} d(n,4,1) = x \log x + (2\gamma - 1 + \log 2 + \pi/2)x + O(x^{1/2})$, (91)

(ii) $\displaystyle\sum_{n \le x} d(n,4,3) = x \log x + (2\gamma - 1 - \log 2 - \pi/2)x + O(x^{1/2})$ (92)

**Theorem 9.2.** ([NK])    Let $x \in \mathbb{R}$, be a large number, and let $a$, $q \in \mathbb{N}$ be fixed integers, $\gcd(a, q) = 1$.

(i) If $.0372831 < a/q < .67181895$, then,

$$\Delta(x,q,a) = \Omega_{-}((x \log x)^{1/4}(\log\log\log x)^{-1/4}).$$ (92)

(ii) If $a/q < .0372831$  or  $a/q > .67181895$, then,

$$\Delta(x,q,a) = \Omega_{+}((x \log x)^{1/4}(\log\log\log x)^{-1/4}).$$ (93)

## 9.3. Other Formula for the Divisor Function Over Arithmetic Progressions

The values of $d(n)$ for the integers $n$ from the residue classes $a \bmod q$, $\gcd(a, q) = 1$, are expected to be equidistributed. This information is expressed by the presence of an small error

$$E(x,q,a) = \sum_{n \le x,\, n \equiv a \bmod q} d(n) - \frac{1}{\varphi(q)} \sum_{n \le x,\, \gcd(n,q)=1} d(n).$$ (94)

An investigation of the divisor function over arithmetic progressions and short interval is carried out in [BS]. In this analysis, the average order of the divisor function is presented. This is attributed to independent works of Selberg and Hooley.

**Theorem 9.3.** ([BS])    For relatively prime integers $a$, $q \in \mathbb{N}$, and for every small number $\varepsilon > 0$, there exists $\delta > 0$ such that the following estimate holds:

$$\sum_{n \le x,\, n \equiv a \bmod q} d(n) = x \frac{P_q(\log x)}{\varphi(q)} + O(\frac{x^{1-\delta}}{\varphi(q)}),$$ (95)

provided that $q < x^{2/3-\varepsilon}$, and the term $P_q(\log x)$  is the residue

$$P_q(\log x) = \underset{s \to 1}{\mathrm{Re}\, s} \left( \zeta^2(s) \frac{x^s}{s} \prod_{p \mid q} (1 - p^{-s}) \right)$$

$$= \frac{\varphi(q)^2}{q^2} \left( \log x + 2\gamma - 1 \right) + \frac{2\varphi(q)}{q} \sum_{p \mid q} \frac{\mu(d) \log d}{d}. \tag{96}$$

## 10. The Harmonic Sum Over Arithmetic Progressions

The generalized Euler constants are the real numbers

$$\gamma(a, q) = -\frac{\psi(a/q) + \log q}{q} \in \mathbb{R}, \tag{97}$$

where $\psi(z) = \Gamma'(z) / \Gamma(z) - \gamma + \sum_{n \le z-1} n^{-1}$ is the logarithmic derivative of the gamma function, see [LR, Theorem 1]. These arise in the asymptotic formula for the harmonic sum over arithmetic progressions.

**Lemma 10.1.** ([LR]) Let $x \in \mathbb{R}$, be a large number, and let $a, q \in \mathbb{N}$ be integers such that $\gcd(a, q) = 1$. Then

$$\sum_{n \le x,\, n \equiv a \bmod q} \frac{1}{n} = \frac{1}{q} \log x + \gamma(a, q) + O(\frac{1}{x}), \tag{98}$$

holds uniformly in $a$, $q$ and $x$ provided that $a < q \le x$.

## 11 Fractional Sums Over Arithmetic Progressions

**Lemma 11.1.** ([SS]) Let $x \in \mathbb{R}$, be a large number, and let $a, q \in \mathbb{N}$ be integers such that $\gcd(a, q) = 1$. Then

$$\sum_{n \le x,\, n \equiv a \bmod q} \{ x/n \} = \frac{1 - \gamma}{q} x + O((\frac{x}{q})^{1/3}), \tag{99}$$

holds uniformly in $a$, $q$ and $x$ provided that $a < q \le x$.

Similar results for other finite sums of fractional parts are described in [PF].

**Conjecture** **11.2.** Let $x \in \mathbb{R}$, be a large number, and let $a, q \in \mathbb{N}$ be integers such that $\gcd(a, q) = 1$. Then

$$\sum_{n \le x, \, n \equiv a \bmod q} \{ x / n \} = \frac{1 - \gamma}{q} x + O((\frac{x}{q})^{1/4 + \varepsilon}), \tag{100}$$

holds uniformly in $a$, $q$ and $x$ provided that $a < q \le x$.